\documentstyle{amsart}

\newtheorem{theorem}{Theorem}
\newtheorem{lemma}[theorem]{Lemma}
\newtheorem{proposition}[theorem]{Proposition}

\theoremstyle{definition}
\newtheorem{definition}[theorem]{Definition}

\theoremstyle{remark}

\newtheorem{example}[theorem]{Example}
\newtheorem{acknowledgement}{Acknowledgement}

\renewcommand\Re{\operatorname{Re}}
\renewcommand\Im{\operatorname{Im}}
\newcommand\grad{\operatorname{grad}}

\newcommand\GL{\operatorname{GL}} %
\newcommand\SL{\operatorname{SL}} 
\newcommand\SU{\operatorname{SU}} %
\newcommand\U{\operatorname{U}} %
\newcommand\id{\operatorname{id}}

\newcommand\eps{\varepsilon}

\newcommand\co{^{\bold C}} 
\newcommand\sst{^{\mathrm s\mathrm s}} 
\newcommand\qu{/\kern-.7ex/} 
\newcommand\sq{\sqrt{-1}}

\begin{document}

\title{Reductive Group Actions on K\"ahler Manifolds}
\author{Eugene Lerman}
\address{University of Pennsylvania, Department of Mathematics,
Philadelphia, Pennsylvania 19104}
\curraddr{University of California, Department of Mathematics, Santa
Cruz, California 95064-1099}
\email{eugene@@cats.ucsc.edu}
\author{Reyer Sjamaar}
\address{University of Pennsylvania, Department of Mathematics,
Philadelphia, Pennsylvania 19104}
\curraddr{Massachusetts Institute of Technology, Department of
Mathematics, Cambridge, Massachusetts 02139-4307}
\email{sjamaar@@math.mit.edu}
\thanks{This work was partially supported by NSF grant DMS92-03398.}
\thanks{The final version will appear elsewhere.}
\date{October 1992}

\maketitle

\begin{abstract}
We prove that the action of a reductive complex Lie group on a
K\"ahler manifold can be linearized in the neighbourhood of a fixed
point, provided that the restriction of the action to some compact
real form of the group is Hamiltonian with respect to the K\"ahler
form.
\end{abstract}

Much attention has been given to the problem of ``dividing out'' a
space by a group of transformations in such diverse areas of
mathematics as classical mechanics and algebraic geometry. In
classical mechanics one typically deals with a phase space $M$, a
surface $L\subset M$ defined as a level set of a number of integrals
of motion, plus a group $G$ of canonical transformations leaving the
motion on $L$ invariant. The ``reduced'' phase space is then the space
$L/G$. Modern treatments of this construction have been given by Meyer
\cite{me:sy} and Marsden and Weinstein \cite{ma:re}.

Algebraic geometers like to study complex-projective varieties $M$
acted upon by groups of projective transformations that are usually
not compact. Often, though, their groups are reductive, i.e., they are
complexifications $G\co$ of compact real Lie groups $G$. A good notion
of a quotient in this category has been defined by Mumford
\cite{mu:ge}. It consists in choosing a subset $M\sst$ of
``semistable'' points of $M$, where the action is ``good'', and then
forming a quotient $M\sst\qu G\co$. This quotient is smaller than the
naive, set-theoretical quotient $M\sst/G\co$, which is not a Hausdorff
space owing to the fact that the orbits of the non-compact group
$G\co$ are not closed. The space $M\sst\qu G\co$ is basically the set
of {\em closed\/} orbits in $M\sst$.

Complex-projective manifolds are not classical phase spaces, but they
have very special symplectic forms, obtained by restricting the
Fubini-Study form on the ambient projective space.  This raises the
question whether there is any relation between the symplectic and the
algebro-geometric quotients of such spaces.  Guillemin, Sternberg,
Mumford and others have made the important observation that for
projective manifolds the notions of a quotient in symplectic and
algebraic geometry are essentially the same. Kempf and Ness
\cite{ke:le} have proved an analogous result for affine varieties. See
Bott \cite{bo:mo} for an account of this story. On the quantum level
this correspondence leads to integral formulae for multiplicities of
representations; see Guillemin and Sternberg \cite{gu:ge}.

We have sought to generalize these multiplicity formulae to the case
where the reduced phase space $L/G$ is not a manifold, but a singular
space. In the process we had to investigate more closely the
correspondence between quotients in mechanics and in algebraic
geometry, and the properties of actions of reductive complex Lie
groups on K\"ahler manifolds. In particular, we had to address the
issue of the existence of slices for such actions. The problem of
linearizing the action near a fixed point and the construction of
slices has been solved by Luna \cite{lu:sl} for affine varieties and
by Snow \cite{sn:re} for Stein spaces. The difficulty of the problem
lies in the fact that an action of a reductive group $G\co$ is never
proper, unless it is locally free. One therefore faces the challenge
of controlling the behaviour of the action ``at infinity in the
group''.  It seems unlikely to us that this can be done in general
without imposing some restrictions on the space $M$.

In this note we present an alternative proof of a theorem of Koras
\cite{ko:li} to the effect that one can linearize a $G\co$-action at a
fixed point on a K\"ahler manifold $M$. The advantage of our proof is
that it can be generalized to show the existence of slices at certain
points of $M$, namely the totally real points. This will be the
subject of a forthcoming paper \cite{le:ho}. As a foretaste we will at
the end of this paper spell out a little example given by Luna,
showing that slices do not exist at arbitrary points of $M$. We are
also preparing a paper on multiplicity formulae for singular reduced
phase spaces.  The result to be discussed here is the following:

\begin{theorem}[Koras]\label{theorem:linear}
Let $M$ be a K\"ahler manifold and let $G\co$ act holomorphically on
$M$. Assume the action of the compact real form $G$ is Hamiltonian.
Let $m$ be any fixed point of the $G\co$-action.  Then the action of
$G\co$ can be linearized in a neighbourhood of $m$\rom, i.e.\rom,
there exist a $G\co$-invariant open neighbourhood $U$ of $m$ in
$M$\rom, a $G\co$-invariant open neighbourhood $U'$ of the origin $0$
in the tangent space $T_mM$ and a biholomorphic $G\co$-equivariant map
$U\to U'$.
\end{theorem}

We find ourselves unable to understand part of Koras' proof of this
theorem. In particular, we fail to see a justification for his
application of the curve selection lemma. Our proof will follow the
same strategy as Koras', using a momentum map for the action of the
compact group $G$ and an analytic continuation argument based on
Weyl's unitary trick. But at the crucial point, instead of invoking
the curve selection lemma, we will rely on an ``interpolation''
argument, deforming the K\"ahler metric in the neighbourhood of a
fixed point.

Let us now explain the statement of the theorem and introduce some
notation. We denote the complex structure on $M$ by $J$, the K\"ahler
metric by $ds^2$ and the K\"ahler form $\Im ds^2$ by $\omega$. Then
$\Re ds^2=\omega(\cdot,J\cdot)$ is the corresponding Riemannian
metric. We may assume without loss of generality that $ds^2$ is
invariant under the compact group $G$. So the transformations on $M$
defined by $G$ are holomorphic and they are isometries with respect to
the K\"ahler metric. Saying that the action of $G$ is Hamiltonian
means that for all $\xi$ in the Lie algebra $\frak g$ of $G$ the
vector field $\xi_M$ on $M$ induced by $\xi$ is Hamiltonian. In this
case we have a {\em momentum map\/} $\Phi$ from $M$ to the dual $\frak
g^*$ of the Lie algebra of $G$ with the property that
$$d\Phi^\xi=\iota_{\xi_M}\omega$$
for all $\xi$. Here $\Phi^\xi$ is the $\xi$-th component of $\Phi$,
defined by $\Phi^\xi(m)=\bigl(\Phi(m)\bigr)(\xi)$. Because its
components are Hamiltonian functions, a momentum map is uniquely
determined up to additive constants. After averaging with respect to
the given action on $M$ and the coadjoint action on $\frak g^*$ we may
assume that the map $\Phi$ is $G$-equivariant. It is easy to give
sufficient conditions for the existence of a momentum map, e.g., the
first Betti number of $M$ is zero, or the K\"ahler form $\omega$ is
exact. (See e.g.\ \cite{gu:sy,we:le}.) More surprisingly, by a theorem
of Frankel \cite{fr:fi} a momentum map always exists if the action has
at least one fixed point and $M$ is {\em compact}.

If $M$ is  $\bold C^n$ with  the standard Hermitian  structure and the
standard symplectic form $\Omega$, then a momentum map is given by the
formula
\begin{equation}\label{equation:quadratic}
\Phi_{\bold C^n}^\xi(v)=1/2\,\Omega\bigl(\xi_{\bold C^n}(v),v\bigr),
\end{equation}
where $\xi _{\bold C^n}$ denotes the image of $\xi \in\frak g$ in
the Lie algebra $\frak s\frak p({\bold C^n},\Omega)$, and $v\in\bold
C^n$. In this setting the statement of the theorem is of course a
tautology, but we will need (\ref{equation:quadratic}) further on.

The decomposition of the complexified Lie algebra $\frak g\co=\frak
g\otimes\bold C$ into a direct sum $\frak g\co=\frak g\oplus\sq\,\frak
g$ gives rise to the Cartan decomposition $G\co=G\exp\sq\,\frak
g$.  (``Polar coordinates'' on the group.) The restriction of the
exponential map $\exp\colon\sq\,\frak g\to G\co$ is a diffeomorphism
onto its image. This means that every element $g$ of $G\co$ can be
uniquely decomposed into a product $g=k\exp\sq\,\xi$, with $k\in G$
and $\xi\in\frak g$.

Because $G\co$ acts holomorphically on $M$, for any $\xi$ in $\frak g$
the vector field $(\sq\,\xi)_M$ induced by $\sq\,\xi$ is equal to
$J\xi_M$. It follows easily from the definition of a momentum map that
$J\xi_M$ is equal to the gradient vector field (with respect to the
Riemannian metric $\Re ds^2$) of the $\xi$-th component of the
momentum map,
\begin{equation}\label{equation:grad}
(\sq\,\xi)_M=J\xi_M=\grad\Phi^\xi.
\end{equation}
This elementary observation will enable us to gain control over the
behaviour of the action ``at infinity in the group''. For one thing,
it implies that the trajectory $\gamma(t)$ of $\grad\Phi^\xi$ through
a point $x$ in $M$ is given by $\gamma(t)=\exp(\sq\,t\xi) x$, which
does not depend on the choice of the K\"ahler metric and the momentum
map.

\begin{trivlist}\item[\hskip\labelsep{\em Proof of Theorem
\ref{theorem:linear}}.]
The proof is in three steps. In the first step we reduce the statement
of the theorem to a statement about the trajectories of the gradient
vector fields $(\sq\,\xi)_M=\grad\Phi^\xi$. In the second step we
consider the case where the K\"ahler metric $ds^2$ is flat in some
neighbourhood of $m$. The third step consists in showing that an
arbitrary metric $ds^2$ can always be deformed to a metric which is
flat close to $m$ and which is still compatible with all the relevant
data.

{\em Step\/} 1. The tangent space $T_mM$ at $m$ is a Hermitian vector
space, which we shall identify with standard $\bold C^n$. Then the
value of the K\"ahler form $\omega$ at $M$ is the standard symplectic
form $\Omega$ on $\bold C^n$. The tangent action of $G\co$ defines a
linear representation $G\co\to\GL(n,\bold C)$, the restriction of
which to $G$ is a unitary representation $G\to\U(n)$. Let $\phi\colon
O\to\bold C^n$ be a local holomorphic coordinate on $M$ with
$\phi(m)=0$ and $d\phi(m)=\id_{\bold C^n}$, where $O$ is a small
$G$-invariant open neighbourhood of $m$ such that $O'=\phi(O)$ is a
ball about the origin in $\bold C^n$. Then the pullback of the form
$\Omega$ is equal to $\omega$ at the point $m$. After averaging over
$G$ and shrinking $O$ if necessary we may assume that $\phi$ is
$G$-equivariant. Let $\psi\colon O'\to O$ be the inverse of $\phi$.

We claim that if $O$ is sufficiently small $\phi$ and $\psi$ can be
uniquely extended to holomorphic $G\co$-equivariant maps
$\phi\co\colon U\to U'$ and $\psi\co\colon U'\to U$, where $U=G\co O$
and $U'=G\co O'$. If this can be done, it is clear that
$\phi\co\circ\psi\co=\id_{U'}$ and $\psi\co\circ\phi\co=\id_{U}$, so
then the theorem will be proved.

Following Heinzner \cite{he:ge} we shall call an open $G$-invariant
subset $A$ of a $G\co$-space $X$ {\em orbitally convex\/} with respect
to the $G\co$-action if for all $x$ in $U$ and all $\xi$ in $\frak g$
the intersection of the curve $\{\,\exp(\sq\,t\xi)x:t\in\bold R\,\}$
with $A$ is connected.

\begin{proposition}[Cf.\ Heinzner \cite{he:ge}, Koras
\cite{ko:li}]\label{proposition:orbit}
Let $X$ and $Y$ be complex manifolds acted upon by $G\co$. If $A$ is
an orbitally convex open subset of $X$ and $f\colon A\to Y$ is a
$G$-equivariant holomorphic map\rom, then $f$ can be uniquely extended
to a $G\co$-equivariant holomorphic map $f\co\colon G\co A\to Y$.

Consequently, if the image $f(A)$ is open and orbitally convex in $Y$
and $f\colon A\to f(A)$ is biholomorphic\rom, then the extension
$f\co\colon G\co A\to Y$ is biholomorphic onto the open subset
$G\co\cdot f(A)$.
\end{proposition}

\begin{pf}
The only way to extend $f$ equivariantly is by putting
$f\co\bigl(g\exp(\sq\,\xi)x\bigr)=g\exp(\sq\,\xi)f(x)$ for all $x$ in
$A$, $g$ in $G$ and $\xi$ in $\frak g$. We have to check this is
well-defined.

Let $x\in A$ and $\xi\in\frak g$ be such that $\exp(\sq\,\xi)x\in A$.
Then by assumption $\exp(\sq\,t\xi)x\in A$ for all $t$ between 0 and
1. So $f\bigl(\exp(\sq\,t\xi)x\bigr)$ is well-defined for $0\leq t\leq
1$.  Define the curves $\alpha(t)$ and $\beta(t)$ in $Y$ by
$\alpha(t)=f\bigl(\exp(\sq\,t\xi)x\bigr)$ and $\beta(t)=
\exp(\sq\,t\xi)f(x)$ for $0\leq t\leq 1$.  Then $\alpha(t)$ and
$\beta(t)$ are integral curves of the vector fields $f_*(\sq\,\xi)_X$
and $(\sq\,\xi)_Y$ respectively, both with the same initial value
$f(x)$. Now since $f$ is $G$-equivariant we have $f_*\xi_X=\xi_Y$,
and, because $f$ is also holomorphic, $f_*(\sq\,\xi)_X= f_*(J\xi_X)=
Jf_*\xi_X= J\xi_Y= (\sq\,\xi)_Y$. Hence $\alpha(t)=\beta(t)$, in other
words $f\bigl(\exp(\sq\,t\xi)x\bigr)=\exp(\sq\,t\xi)f(x)$ for $0\leq
t\leq 1$.

It follows that for all $x$ in $A$ and all $\xi$ in $\frak g$ such
that $\exp(\sq\,\xi)x$ is in $A$ we have
$f\bigl(\exp(\sq\,\xi)x\bigr)=\exp(\sq\,\xi)f(x)$. It is easy to
deduce from this that $f\co$ is well-defined.

Finally observe that if the image $f(A)$ is open and orbitally convex
in $Y$ and $f\colon A\to f(A)$ is biholomorphic, then the inverse
$f^{-1}$ also has a holomorphic extension $(f^{-1})\co\colon G\co
f(A)\to G\co A$, and by uniqueness this must be the inverse of $f\co$.
\end{pf}

This proposition reduces the proof of Theorem \ref{theorem:linear} to
showing that $O\subset M$ and $O'\subset\bold C^n$ are orbitally
convex.

{\em Step\/} 2. We will now show that $O$ is orbitally convex if the
K\"ahler metric is flat on an open subset containing $O$. Since
$\bold C^n$ is flat, this covers the case of $O'\subset\bold C^n$, so
in particular we will have shown $O'$ is orbitally convex.

The argument is a variation on a convexity argument of Kempf and Ness
\cite{ke:le}. The assumption that the K\"ahler metric is flat on an
open subset containing $O$ implies that the holomorphic coordinate map
$\phi\colon O\to\bold C^n$ is a K\"ahler isometry onto the open ball
$O'$ about the origin. After choosing appropriate constants we may
assume that $\Phi(m)=0$. Then $\psi^*\Phi\colon O'\to\frak g$ is equal
to the quadratic momentum map given by (\ref{equation:quadratic}).
(Recall $\psi=\phi^{-1}$.) In order not to overburden the notation we
will in the proof of the following lemma identify $O$ with $O'$ and
$\Phi$ with the momentum map (\ref{equation:quadratic}). By $R(x)$ we
denote the Riemannian distance of $x\in O$ to $m$.

\begin{lemma}\label{lemma:angle}
Assume the K\"ahler metric is flat in the neighbourhood of $O$. Then
for all $\xi\in\frak g$ and $v\in O$ the momentum function $\Phi^\xi$
measures the inner product of the outward pointing normal $\grad R^2$
to the metric sphere of radius $R$ about $m$ and the vector field
$J\xi_O=\grad\Phi^\xi$\rom, as follows\rom:
\begin{equation}\label{equation:angle}
\bigl\langle\grad R^2(v),\grad\Phi^\xi(v)\bigr\rangle=4\Phi^\xi(v).
\end{equation}

It follows that $O$ is orbitally convex with respect to the
$G\co$-action.
\end{lemma}

\begin{pf}
Let $\delta(t)=\exp(\sq\,t\xi)v$ denote the gradient trajectory of
$\Phi^\xi$ through $v$. On one hand,
$$
\frac{d}{dt}R^2\bigl(\delta(t)\bigr) = \bigl\langle\grad
R^2\bigl(\delta(t)\bigr),\delta'(t)\bigr\rangle = \bigl\langle\grad
R^2\bigl(\delta(t)\bigr),\grad\Phi^\xi\bigl(\delta(t)\bigr)\bigr\rangle.
$$
On the other hand,
\begin{align*}
\frac{d}{dt}R^2\bigl(\delta(t)\bigr)&=
\frac{d}{dt}\bigl\|\delta(t)\bigr\|^2
= \frac{d}{dt}\bigl\langle\delta(t),\delta(t)\bigr\rangle =
2\bigl\langle\delta'(t),\delta(t)\bigr\rangle = \\
&=
2\bigl\langle(\sq\,\xi)_O\bigl(\delta(t)\bigr),\delta(t)\bigr\rangle =
2\Omega\bigl(\xi_O\bigl(\delta(t)\bigr),\delta(t)\bigr) =
4\Phi^\xi\bigl(\delta(t)\bigr),
\end{align*}
where we have used (\ref{equation:quadratic}) and
(\ref{equation:grad}). Taking $t=0$ yields (\ref{equation:angle}).

Now (\ref{equation:angle}) implies that the curve $\delta(t)$ can only
enter $O$ at a point $p$ in the boundary $\partial O$ for which
$\Phi^\xi(p)\leq0$ and leave it at a point $q\in\partial O$ where
$\Phi^\xi(q)\geq0$. But $\delta(t)$ is also a gradient curve of the
function $\Phi^\xi$ and so $\Phi^\xi$ is increasing along $\delta(t)$
for all $t\in\bold R$. If $\delta(t)$ is not constant,
$\Phi^\xi\bigl(\delta(t)\bigr)$ is strictly increasing. Therefore, if
$\delta(t)$ leaves the ball $O$ at some point, it can never sneak back
in. Consequently $\{\,\delta(t):t\in\bold R\,\}\cap O$ is connected.
If $\delta(t)$ is constant it is trivially true that
$\{\,\delta(t):t\in\bold R\,\}\cap O$ is connected.
\end{pf}

{\em Step\/} 3. We claim there exists a K\"ahler metric $d\tilde s^2$
on $M$ with the following properties: It is compatible with the
holomorphic structure $J$ and $G$-equivariant, the $G$-action is
Hamiltonian with respect to $d\tilde s^2$, and $d\tilde s^2$ is flat
in a neighbourhood of $m$. Together with the result of step 2 this
will conclude the proof of Theorem \ref{theorem:linear}.

Let $u$ be a K\"ahler potential for $ds^2$ in a neighbourhood of
$m$, i.e.,
$$
\omega=\sq\,\partial\bar\partial u.
$$
Consider the Taylor expansion $u(x)\sim\sum_{k,l}c_{kl}x^k\bar x^l$ of
$u$ about $m$, where $x$ is the system of holomorphic coordinates
defined by the map $\phi\colon O\to\bold C^n$, and where $k$ and $l$
denote multi-indices. The value of $\omega$ at $m$ is the standard
symplectic form on $\bold C^n$, so $\partial^2u(m)/\partial
x_\alpha\partial\bar x_\beta=\delta_{\alpha\beta}$. Also, deleting the
holomorphic terms $\sum_{k}c_{k0}x^k$ and the antiholomorphic terms
$\sum_{l}c_{0l}\bar x^l$ up to any finite order does not affect
$\partial\bar\partial u$, so we may assume
$$
u(x)=\|x\|^2+\cal R(x),
$$
with remainder term $\cal R(x)$ of order $O\bigl(\|x\|^3\bigr)$. Let
$\chi\colon\bold R\to[0,1]$ be a smooth function with $\chi(t)=0$ for
$t\leq1$ and $\chi(t)=1$ for $t\geq2$. Pick $\lambda>0$ and put
$$
\tilde u(x)=\|x\|^2+ \chi\bigl(\|x\|^2/\lambda^2\bigr)\cal R(x).
$$
Let $O(r)$ be the ball $\{\,x\in O:\|x\|^2<r\,\}$ and define a
smooth two-form $\tilde\omega$ on $M$ by
$$
\tilde\omega=\left\{
\begin{array}{ll}
\sq\,\partial\bar\partial\tilde u&\mbox{ on }
O(3\lambda^2) \\
\omega&\mbox{ on }M-O(2\lambda^2).
\end{array}\right.
$$
Then on $O(\lambda^2)$ the form $\tilde\omega$ is equal to the
standard symplectic form $\Omega$ and so is $G$-invariant there. On
$O(3\lambda^2)$ we have $\tilde\omega-\omega=
\sum_{\alpha,\beta=1}^nf_{\alpha\beta}\,dx_\alpha\wedge d\bar x_\beta$
with
$$
f_{\alpha\beta}(x) = \frac{\partial^2}{\partial
x_\alpha\partial\bar
x_\beta}\bigl(\chi\bigl(\|x\|^2/\lambda^2\bigr)-1\bigr)\cal R(x).
$$
By carrying out the differentiation and using $\chi(t)=1$ for $t\geq2$
and $\cal R(x)=O\bigl(\|x\|^3\bigr)$ one can check in a
straightforward manner that for every $\eps>0$ there exists
$\Lambda>0$ such that for all $\lambda<\Lambda$ and for all $x\in
O(3\lambda^2)$ one has $|f_{\alpha\beta}(x)|<\eps$. In other words,
$\tilde\omega$ can be made arbitrarily close to $\omega$ uniformly on
$M$.  It follows that for $\eps$ and $\lambda$ small enough the
symmetric bilinear form $\tilde\omega(\cdot,J\cdot)$ is
positive-definite, and therefore $\tilde\omega$ is the imaginary part
of a K\"ahler metric $d\tilde s^2$.  By construction $d\tilde s^2$ is
equal to $ds^2$ on $M-O(2\lambda^2)$ and flat on $O(\lambda^2)$. After
averaging over $G$ we may assume $d\tilde s^2$ is $G$-invariant
everywhere. (This does not destroy the flatness on $O(\lambda^2)$,
because $d\tilde s^2$ was already invariant there.)

It remains to be proved that the $G$-action is Hamiltonian with
respect to $\tilde\omega$. By the Poincar\'e Lemma we can find an
$\tilde\omega$-momentum map $\tilde\Phi\colon O(3\lambda^2)\to\frak
g^*$. On the set $O(3\lambda^2)-O(2\lambda^2)$ we have
$\tilde\omega=\omega$, so there $\tilde\Phi$ differs by a locally
constant function $c$ from the $\omega$-momentum map $\Phi$. Since
$O(3\lambda^2)-O(2\lambda^2)$ is connected, $c$ is a constant.
Shifting $\Phi$ by this constant we can paste $\Phi$ and $\tilde\Phi$
together to obtain a global $\tilde\omega$-momentum map for the
$G$-action.\qed
\end{trivlist}

In \cite{le:ho} we will show how to construct slices for a
$G\co$-action on a K\"ahler manifold at points $m$ for which the real
orbit $Gm$ is {\em totally real}. (This includes the case of a fixed
point). Let us recall the definition of a slice.

\begin{definition}
A {\em slice\/} at $m$ for the $G\co$-action is a locally closed
analytic subspace $S$ of $M$ with the following properties:
\begin{enumerate}
\item $m\in S$;
\item the saturation $G\co S$ of $S$ is open in $M$;
\item $S$ is invariant under the action of the stabilizer $(G\co)_m$;
\item\label{bundle} the natural $G\co$-equivariant map from the
associated bundle $G\co\times_{(G\co)_m}S$ into $M$, which sends an
equivalence class $[g,y]$ to the point $gy$, is an analytic
isomorphism onto $G\co S$.
\end{enumerate}
\end{definition}

It follows from (\ref{bundle}) that if a slice at $m$ exists, then for
all $y$ in a $G\co$-invariant neighbourhood of $m$ the stabilizer
$(G\co)_y$ is conjugate to a subgroup of $(G\co)_m$. Furthermore, a
slice has to be nonsingular at $m$ and transverse to the orbit $G\co
m$. Theorem \ref{theorem:linear} clearly shows that slices exist at
fixed points, where $(G\co)_m$ is equal to $G\co$.

The basic idea of the proof of the slice theorem in \cite{le:ho} is
the same as that of Theorem \ref{theorem:linear}: We show that a
totally real $G$-orbit has sufficiently small orbitally convex
neighbourhoods and then use analytic continuation. The main technical
difficulty is to make the interpolation argument of step 3 work
globally along a totally real orbit. This involves a closer scrutiny
of the second-order operator $\partial\bar\partial$.

The following example shows that slices need not exist everywhere,
even at points $m$ where the action is free, that is, $(G\co)_m=1$.

\begin{example}[Luna \cite{lu:sl}]
Let $V$ be the vector space of homogeneous cubic polynomials in two
indeterminates $x$ and $y$. Let $G$ be the group $\SU(2)$; then
$G\co=\SL(2,\bold C)$. Both groups act by linear substitutions in the
two variables. If we declare the monomials $x^3$, $x^2y$, $xy^2$ and
$y^3$ to be an orthonormal basis, we obtain a $G$-invariant Hermitian
inner product on $V$. Every polynomial in $V$ can be written as a
product of three linear factors, which are unique up to scalar
multiples. So, obviously, every element of $G\co$ fixing the
polynomial has to permute these factors (up to scalars). We can
distinguish between four types of polynomials according to the number
of distinct factors in the factorization, and it is easy to compute
their stabilizers. (See table \ref{table:stabilizer}.) The points of
type I form clearly a dense subset of $V$.
\begin{table}\caption{}\label{table:stabilizer}
\begin{center}
\begin{tabular}{|r|c|c|c|}
\hline
&{\em number of distinct factors}&$m$&$(G\co)_m$\\
\hline\hline
I&3&$(x+y)(x-y)y$&$\bold Z/3$\\ \hline
II&2&$x^2y$&1\\ \hline
III&1&$y^3$&$\bold Z/3$-extension of $\bold C$\\ \hline
IV&0&0&$G\co$\\ \hline
\end{tabular}
\end{center}
\end{table}
Consider the point $m=x^2y$. The action is free at $m$, so the orbit
$G\co m$ is a copy of the three-dimensional $\SL(2,\bold C)$ embedded
in the four-dimensional $V$. A straightforward computation shows the
orthogonal complement to $G\co m$ at $m$ is the one-dimensional
subspace generated by $y^3$. We can define a holomorphic
$G\co$-equivariant map $\phi$ from the associated bundle
$E=G\co\times_{(G\co)_m}\bold C=\SL(2,\bold C)\times\bold C$ into $V$
by putting $\phi\bigl([g,\eta]\bigr)=g\cdot(x^2y+\eta y^3)$. The
differential of $\phi$ is bijective at all points of the zero section
of $E$ and the restriction of $\phi$ to the zero section is an
embedding. Therefore $\phi$ is an open immersion in a neighbourhood of
the zero section. Does there exist a $G\co$-invariant neighbourhood
$U$ of the zero section of $E$ such that the restriction of $\phi$ to
$U$ is biholomorphic onto its image? We could draw this conclusion if
we knew that $\phi$ was injective on a $G\co$-invariant neighbourhood
of the zero section. But there are points of type I arbitrarily close
to $m$, whose stabilizer consists of three elements and so cannot be
conjugate to a subgroup of $(G\co)_m=1$. Hence there cannot possibly
exist a slice at $m$, and therefore the map $\phi$ cannot be injective
in any $G\co$-invariant neighbourhood of the zero section of $E$. The
reason is that $\phi$ plays havoc with the fibres of $E$. Let us write
$\eta=-\eps^2$; then $\phi\bigl([g,\eta]\bigr)= gm_\eps$ with
$m_\eps=(x+\eps y)(x-\eps y)y$. Let $U$ be any $\SL(2,\bold
C)$-invariant neighbourhood of the zero section in $E$.  Then $U$
contains an invariant open neighbourhood of the form
$E_\Delta=\SL(2,\bold C)\times\Delta$, where $\Delta$ is a small disc
about the origin in $\bold C$. Choose $\eps$ such that
$\eps^2\in\Delta$. The stabilizer $(G\co)_{m_\eps}$ is generated by
the following matrix of order three:
$$
A_\eps=
\begin{pmatrix}
-1/2&3\eps/2\\-1/2\eps&-1/2
\end{pmatrix},
$$
and we have
$$
A_\eps m= -\frac{1}{8\eps}x^3+\frac{5}{8}x^2y-
\frac{3\eps}{8}xy^2-\frac{9\eps^2}{8}y^3.
$$
So $A_\eps$ moves $m$ way out along its orbit while leaving $m_\eps$
fixed. We conclude that the images of the fibres
$\{I\}\times\Delta\subset E_\Delta$ and $\{A_\eps\}\times\Delta\subset
E_\Delta$ under the map $\phi$ intersect at the point $m_\eps$.
\end{example}

\begin{acknowledgement}
We are grateful to Eugenio Calabi for showing us how to glue K\"ahler
potentials.
\end{acknowledgement}

\makeatletter \renewcommand{\@biblabel}[1]{\hfill#1.}\makeatother

\end{document}